\documentclass[12pt]{article}
\usepackage{mystyle}
\begin{document}
\title{The internal relation}
\author{Gabriel Goldberg}
\maketitle
\begin{abstract}
This paper explores various generalizations of the Mitchell order, focusing mostly on a generalization called the internal relation. The internal relation lacks the implicit strength requirement in the definition of the Mitchell order, and therefore can fail to be wellfounded. We establish some constraints on the illfoundedness of the internal relation, which leads to a proof of a conjecture of Steel regarding rank-to-rank cardinals.
\end{abstract}
\section{Introduction}
This paper explores various generalizations of the Mitchell order, focusing mostly on a generalization called the internal relation. This relation turns out to be related to a number of questions in large cardinal theory, so let us just give some examples. We prove converses to the commuting ultrapower lemma of Kunen. We show that a nonprincipal countably complete ultrafilter on a cardinal can be amenable to its own ultrapower. We prove a conjecture of Steel \cite{SteelWellfounded} regarding the Mitchell order on rank-to-rank extenders. We show that if \(U\) and \(W\) are normal ultrafilters on \(\kappa\) with \(U\mo W\), then \(j_U(\alpha)\leq j_W(\alpha)\) for all ordinals \(\alpha\). We prove that if \(j_0,j_1 : M\to N\) are elementary embeddings between two inner models and \(j_0\) is definable over \(M\), then \(j_0(\alpha) \leq j_1(\alpha)\) for all ordinals \(\alpha\). Finally we analyze the internal relation on countably complete ultrafilters assuming a principle called the Ultrapower Axiom. This analysis is important in \cite{Frechet} and \cite{SC}.

\section{The internal relation}
We now define the internal relation. 
\begin{defn}\label{InternalDefn}
The {\it internal relation \(\I\)} is defined on extenders \(E\) and \(F\) by setting \(E\I F\) if \(j_E\restriction x\in M_F\) for all sets \(x\in M_F\). 
\end{defn}

The definition of the internal relation should be contrasted with that of the generalized Mitchell order:
\begin{defn}
The {\it Mitchell order \(\mo\)} is defined on extenders \(E\) and \(F\) by setting \(E\mo F\) if \(E\in M_F\).
\end{defn}
If \(U\) and \(W\) are countably complete ultrafilters on cardinals \(\lambda_U\) and \(\lambda_W\), then for \(U\mo W\) to hold, \(W\) must have a certain amount of strength: \(P(\lambda_U)\) must belong to \(M_W\). On the other hand, there is no such strength requirement for the internal relation. We will see examples below (for example in \cref{CommutingUltrapowers}) of countably complete uniform ultrafilters \(U\I W\) where \(P(\lambda_U)\notin M_W\) so \(U\not \mo W\).

It is also not the case that \(E\mo F\) implies \(E \I F\). The problem is that even if \(E\in M_F\), it need not be true that \(j_E^{M_F} = j_E\restriction M_F\). (It is even possible to have \(E\mo F\)  and \(E\I F\) but \(j_E^{M_F} \neq j_E\restriction M_F\).)

For these reasons the behavior of the internal relation is quite different from that of the Mitchell order. For example, we will see in \cref{CommutingUltrapowers} that while the internal relation is irreflexive on nontrivial extenders, it actually has 2-cycles.

\subsection{Locality of the internal relation}
In \cref{InternalDefn}, we set \(E\I F\) if \(j_E\restriction M_F\) was an amenable class of \(M_F\). We now show that this implies that in fact \(j_E\restriction M_F\) is a definable class of \(M_F\), and indeed \(j_E\restriction M_F\) is an extender embedding of \(M_F\). Thus the internal relation is {\it local} in the sense that it depends only on the existence of a certain set in \(M_F\) (namely the extender giving rise to \(j_E\restriction M_F\)), and not on the amenability of an entire class. 
\begin{lma}\label{sExtender}
Suppose \(E\) and \(F\) are extenders and \(E\I F\). Then there is an extender \(E'\) of \(M_F\) such that \(j^{M_F}_{E'} = j_E\restriction M_F\).
\begin{proof}
It suffices to show that \(j_E(M_F) = H^{j_E(M_F)}(j_E[M_F]\cup \lambda)\) for some ordinal \(\lambda\). Since \(F\) is an extender, there is an ordinal \(\lambda_F\) such that \(M_F = H^{M_F}(j_F[V]\cup \lambda_F)\). By elementarity, \begin{equation}\label{Hull1} j_E(M_F) = H^{j_E(M_F)}(j_E(j_F)[M_E]\cup j_E(\lambda_F))\end{equation} Since \(E\) is an extender, there is an ordinal \(\lambda_E\) such that \(M_E = H^{M_E}(j_E[V]\cup \lambda_E)\). Thus \begin{align}j_E(j_F)[M_E] &= H^{j_E(M_F)}(j_E(j_F)[j_E[V]]\cup j_E(j_F)[\lambda_E])\nonumber\\ &= H^{j_E(M_F)}(j_E\circ j_F[V]\cup j_E(j_F)[\lambda_E])\label{Hull2}\end{align}
Combining \cref{Hull1} and \cref{Hull2}, we obtain 
\begin{align*}j_E(M_F) &= H^{j_E(M_F)}(j_E\circ j_F[V]\cup j_E(j_F)[\lambda_E]\cup j_E(\lambda_F))\\
&\subseteq H^{j_E(M_F)}(j_E[M_F]\cup j_E(j_F)(\lambda_E)\cup j_E(\lambda_F))
\end{align*}
Letting \(\lambda = \sup  \{j_E(j_F)(\lambda_E), j_E(\lambda_F)\}\), we have \(j_E(M_F) = H^{j_E(M_F)}(j_E[M_F]\cup \lambda)\), as desired.
\end{proof}
\end{lma}

We take this opportunity to set up some notation which allows us to point out one of the main distinctions between ultrafilters and extenders in this context.
\begin{defn}
Suppose \(M_0, M_1,\) and \(N\) are models of ZFC. We write \[(i_0,i_1) : (M_0,M_1)\to N\] to denote that \(i_0 : M_0\to N\) and \(i_1 : M_1\to N\) are elementary embeddings.
\end{defn}

\begin{defn}
We say \((i_0,i_1) : (M_0,M_1)\to N\) is {\it minimal} if \(N = H^N(i_0[M_0]\cup i_1[M_1])\).
\end{defn}

It is not hard to see that if \((i'_0,i'_1) : (M_0,M_1)\to N'\), there is a unique minimal \((i_0,i_1) : (M_0,M_1)\to N\) such that there is an elementary \(h : N\to N'\) with \(i_0' = h\circ i_0\) and \(i_1' = h \circ i_1\).

The following lemma follows from the proof of \cref{sExtender}, but the lemma is just false for extenders in general.

\begin{lma}\label{InternalMinimal}
Suppose \(U\) and \(W\) are ultrafilters. Then \[(j_{U}(j_{W}),j_{U}\restriction M_{W}) : (M_{U},M_{W})\to j_{U}(M_{W})\] is minimal.
\begin{proof}
Let \(a = [\text{id}]_{W}\). Then \(M_{W} = H^{M_{W}}(j_{W}[V]\cup \{a\})\). Therefore by the elementarity of \(j_U\), \begin{align*}j_{U}(M_{W}) &= H^{j_{U}(M_{W})}(j_{U}(j_{W})[M_{U}]\cup \{j_{U}(a)\})\\& \subseteq  H^{j_{U}(M_{W})}(j_{U}(j_{W})[M_{U}]\cup j_{U}[M_{W}])\end{align*}
In other words, \((j_{U}(j_{W}),j_{U}\restriction M_{W}) : (M_{U},M_{W})\to j_{U}(M_{W})\) is minimal.
\end{proof}
\end{lma}

In \cite{RF} Lemma 5.5, we show that if \(j_0 : V\to M_0\) and \(j_1 : V\to M_1\) are ultrapower embeddings, \((i_0,i_1) : (M_0,M_1)\to N\) is minimal, and \(i_0 \circ j_0 = i_1\circ j_1\), then \(i_0\) and \(i_1\) are ultrapower embeddings of \(M_0\) and \(M_1\) respectively. It follows that if \(U\I W\), then \(j_{U}\restriction M_{W}\) is an ultrapower embedding of \(M_{W}\). We now study a particular \(M_{W}\)-ultrafilter that gives rise to \(j_{U}\restriction M_{W}\).

\begin{defn}
An ultrafilter \(U\) on an ordinal \(\alpha\) is {\it tail uniform} (or just {\it uniform}) if \(\alpha\setminus \beta\in U\) for all \(\beta < \alpha\). If \(U\) is a uniform ultrafilter, then \(\textsc{sp}(U)\) denotes the underlying set of \(U\), i.e., the unique ordinal \(\alpha\) that belongs to \(U\). 
\end{defn}

\begin{defn}
Suppose \(U\) and \(W\) are countably complete ultrafilters with \(U\) uniform on an ordinal \(\lambda\). Then \(\s W U = \{X\in P^{M_W}(\lambda_*) : j_W^{-1}[X]\in U\}\) where \(\lambda_* = \sup j_W[\lambda]\).
\end{defn}

\begin{lma}\label{Pushforward}
Suppose \(U\) and \(W\) are countably complete ultrafilters with \(U\) uniform on an ordinal \(\lambda\). Then \(\s W U\) is the uniform \(M_W\)-ultrafilter derived from \(j_U\restriction M_W\) using the ordinal \(j_U(j_W)([\textnormal{id}]_U)\). Moreover \(j_{\s W U}^{M_W} = j_U\restriction M_W\).
\begin{proof}
This is a simple calculation: for any \(X\subseteq \lambda_*\) with \(X\in M_W\),
\begin{align*}
X\in \s W U &\iff j_W^{-1}[X]\in U\\
&\iff \{\alpha < \lambda : j_W(\alpha)\in X\}\in U\\
&\iff j_U(j_W)([\text{id}]_U)\in j_U(X)
\end{align*}
The final equivalence follows from Los's theorem.

The fact that \(j_{\s W U}^{M_W} = j_U\restriction M_W\) follows from the fact that \[j_U(M_W) = H^{j_U(M_W)}(j_U[M_W]\cup j_U(j_W)([\text{id}]_U))\] This is an immediate consequence of the minimality of \[(j_U(j_W),j_U\restriction M_W) : (M_U,M_W)\to j_U(M_W)\] proved in \cref{InternalMinimal}.
\end{proof}
\end{lma}
\subsection{An ultrafilter amenable to its own ultrapower}
In this section we briefly study another generalization of the Mitchell order that does away with the implicit strength requirement of the Mitchell order in the most dramatic way possible.
\begin{defn}
The {\it amenability relation \(\A\)} is defined on countably complete ultrafilters \(U\) and \(W\) by setting \(U\A W\) if \(U\cap M_W\in M_W\).
\end{defn}

Note that the amenability relation is not invariant under isomorphisms. For example, if no set \(X\in U\) belongs to \(M_W\), then \(U\A W\) simply because \(U\cap M_W = \emptyset\). In fact this sensitivity makes all the difference between the amenability relation and the internal relation:

\begin{lma}
Suppose \(U\) and \(W\) are countably complete ultrafilters. Then \(U\I W\) if and only if \(U'\A W\) for all ultrafilters \(U'\) isomorphic to \(U\).
\begin{proof}
Suppose \(U\I W\). Suppose \(U'\) is isomorphic to \(U\). We will show \(U'\A W\). We may assume without loss of generality that for some \(X\in U'\), \(X\in M_W\), since otherwise \(U'\A W\) trivially. Thus \([\text{id}]_{U'} \in j_{U'}(X)\subseteq j_{U'}(M_W) = j_{U}(M_W)\subseteq M_W\). Since \(U'\cap M_W\) is the ultrafilter derived from \(j_U\restriction M_W\) using \([\text{id}]_{U'}\), we can define \(U'\cap M_W\) inside \(M_W\) using the definable class  \(j_U\restriction M_W\)  and the parameter  \([\text{id}]_{U'}\). Therefore \(U'\cap M_W\in M_W\), so \(U'\A W\) as desired.

Conversely assume \(U'\A W\) for all \(U'\) isomorphic to \(U\). We must show \(U\I W\). Without loss of generality we may assume \(U\) is a uniform ultrafilter on an ordinal \(\lambda\). Let \(U' = (j_W)_*(U) = \{X\subseteq \sup j_W[\lambda] : j_W^{-1}[X]\in U\}\), the pushforward of \(U\) via the function \(j_W: \lambda \to \sup j_W[\lambda]\). Then \(U'\) is isomorphic to \(U\) and so \(U'\cap M_W\in M_W\). But \(U'\cap M_W = \s W U\). It follows from \cref{Pushforward} that \(U\I W\).
\end{proof}
\end{lma}

Given the well-known fact that for any nonprincipal countably complete ultrafilter \(U\), \(U\notin M_U\), it is natural to ask whether it is possible that \(U\cap M_U\in M_U\). One might naively expect that the amenability relation is irreflexive on nonprincipal countably complete ultrafilters. Here we will show that this is false assuming the existence of a supercompact cardinal, and even a bit less.

\begin{thm}\label{SelfAmenable}
Suppose \(\kappa\) is a supercompact cardinal. Then there is a nonprincipal \(\kappa\)-complete ultrafilter \(U\) such that \(U\cap M_U\in M_U\).
\begin{proof}
Fix \(\lambda\) such that \(\text{cf}(\lambda)\geq \kappa\) and \(2^\lambda = \lambda^+\). (By Solovay's theorem, any singular strong limit cardinal \(\lambda\) of cofinality at least \(\kappa\) will suffice.) Fix a \(\kappa\)-complete weakly normal ultrafilter \(D\) on \(\lambda\) such that \(\lambda^+\) carries a uniform \(\kappa\)-complete ultrafilter in \(M_D\). In \(M_D\), let \(Z\) be a weakly normal, \(j_D(\kappa)\)-complete, \(j_D(\lambda^+)\)-supercompact ultrafilter on \(j_D(\lambda^+)\) such that letting \(N = M_Z^{M_D}\), \(o^N(j_D(\lambda^+)) = [\text{id}]_D\). Let \(i : V\to N\) be the composition \(i = j_Z^{M_D}\circ j_D\).
\begin{clm}
\(N = H^N(i[V]\cup \{\sup i[\lambda]\})\).
\begin{proof}
Note that \(\sup i [\lambda^+] = \sup j_Z^{M_D}[j_D(\lambda^+)]\). Therefore 
\begin{align*}N &= H^N(j_Z^{M_D}[M_D]\cup \{[\text{id}]_Z\}) \\&= H^N(i[V]\cup \{\sup i[\lambda], j_Z^{M_D}([\text{id}]_D)\})\end{align*}

By Solovay's Lemma, \(j_Z^{M_D}[j_D(\lambda^+)]\) is definable in \(N\) from \(\sup i [\lambda^+]\) and \(j_Z^{M_D}(\vec T)\) where \(\vec T\) is any stationary partition of the set of cofinality \(\omega\) ordinals below \(j_D(\lambda^+)\) into \(j_D(\lambda^+)\)-many pieces. In particular if  \(\vec S\) is a stationary partition of the set of cofinality \(\omega\) ordinals below \(\lambda^+\) into \(\lambda^+\)-many pieces then \(j_Z^{M_D}[j_D(\lambda^+)]\) is definable in \(N\) from \(\sup i [\lambda^+]\) and \(j_Z^{M_D}(j_D(\vec S)) = i(\vec S)\). So \(j_Z^{M_D}[j_D(\lambda^+)]\in H^N(i[V]\cup \{\sup i[\lambda]\})\).

Note also that \(j_D(\lambda^+) = \text{cf}^N(\sup i [\lambda^+]) \in H^N(i[V]\cup \{\sup i[\lambda]\})\). Since \(o^N(j_D(\lambda^+)) = [\text{id}]_D\), it follows that \([\text{id}]_D\in H^N(i[V]\cup \{\sup i[\lambda]\})\).

Since \([\text{id}]_D\) and \(j_Z^{M_D}[j_D(\lambda^+)]\) both belong to \(H^N(i[V]\cup \{\sup i[\lambda]\})\), so does \(j_Z^{M_D}([\text{id}]_D)\). Therefore \[N = H^N(i[V]\cup \{\sup i[\lambda], j_Z^{M_D}([\text{id}]_D)\}) = H^N(i[V]\cup \{\sup i[\lambda]\})\] as claimed.
\end{proof}
\end{clm}

The following claim is standard, due perhaps to Menas \cite{Menas}, and we omit the proof:

\begin{clm}\label{Covering}
Every set \(A\subseteq N\) such that \(|A| \leq \lambda^+\) is covered by a set \(A'\in N\) such that \(|A'| < i(\kappa)\).\qed
\end{clm}

Let \(W\) be the ultrafilter derived from \(i\) using \(\sup i[\lambda]\). Then \(W\) is weakly normal, \(M_W = N\), and \(j_W =  i\).

\begin{clm}\label{Extension}
Suppose \(F\) is a uniform \(\kappa\)-complete filter on \(\lambda^+\) that is generated by some \(H\subseteq F\) with \(|H| = \lambda^+\). Then \(F\) extends to an ultrafilter \(U\) that is isomorphic to \(W\).
\begin{proof}
Again this is a standard fact using \cref{Covering}. Note that \(i[H]\) is covered by a set \(H'\in N\) such that \(|H'| < i(\kappa)\). Since \(i(F)\) is \(i(\kappa)\)-complete, \(\bigcap H'\cap i(F)\) is nonempty. Fix \(\xi\in \bigcap H'\cap i(F)\). Note that \(\xi \in \bigcap i[H]\). Therefore the ultrafilter \(U\) derived from \(i\) using \(\xi\) extends \(H\). Moreover \(U\) is a uniform ultrafilter on \(\lambda^+\) and \(U\RK W\). Since \(W\) is a weakly normal ultrafilter on \(\lambda^+\), \(W\) is minimal in the Rudin-Keisler order on uniform ultrafilters on \(\lambda^+\). It follows that \(U\) and \(W\) are isomorphic.
\end{proof} 
\end{clm}

Let \(H\) be a uniform \(\kappa\)-complete ultrafilter on \(\lambda^+\) in \(M_D\) and let \(F\) be the filter generated by \(H\) in \(V\). Then \(F\) is a \(\kappa\)-complete uniform filter on \(\lambda^+\). Note that \(|P(\lambda^+)\cap M_D| \leq |j_D(\kappa)| \leq \kappa^\lambda = 2^\lambda = \lambda^+\). Therefore \(|H| = \lambda^+\). Applying \cref{Extension}, let \(U\) be an extension of \(F\) that is isomorphic to \(W\). Then \[U \cap M_U = U\cap M_W = F\in M_W = M_U\] So \(U\) is as desired.
\end{proof}
\end{thm}

Note that the proof of the theorem shows that ultrafilters witnessing failures of irreflexivity in the amenability relation can be relatively simple:
\begin{cor}[GCH]
Suppose \(\kappa\) is \(\kappa^{++}\)-supercompact. Then there is a \(\kappa\)-complete ultrafilter \(U\) on \(\kappa^+\) such that \(U\cap M_U\in M_U\). Moreover \(U\) is isomorphic to a weakly normal ultrafilter on \(\kappa^+\).\qed
\end{cor}
We do not know whether a weakly normal countably complete ultrafilter can be amenable to its own ultrapower.
\subsection{Cycles in the internal relation}
We briefly prove the standard fact that no extender satisfies \(E\I E\). 
\begin{lma}
The internal relation is irreflexive on nontrivial extenders.
\begin{proof}
Suppose \(E\) is an extender and \(E\I E\). We will show that \(E\) is trivial, or in other words \(j_E\) is the identity. Note that since \(E\I E\), \(j_E\restriction \alpha\) is in \(M_E\) for all ordinals \(\alpha\). It follows that \(P(\alpha)\subseteq M_E\) for all ordinals \(\alpha\), since for any \(X\subseteq \alpha\), \(X =\{\beta < \alpha : (j_E\restriction \alpha) (\beta) \in  j_E(X)\}\). Since every set of ordinals belongs to \(M_E\), \(M_E = V\). Therefore \(j_E\) is the identity by Kunen's inconsistency theorem \cite{Kunen}, so \(E\) is trivial.
\end{proof}
\end{lma}

In this section, we give examples of 2-cycles in the internal relation; i.e., extenders \(E\) and \(F\) such that \(E\I F\) and \(F\I E\). The first examples we discuss come from Kunen's commuting ultrapowers lemma (see \cite{Larson}).
\begin{thm}[Kunen]\label{CommutingUltrapowers}
Suppose \(U\) is a countably complete ultrafilter on a set \(X\) and \(W\) is an ultrafilter that is closed under intersections of \(X\)-sequences. Then \(j_U(j_W) = j_W\restriction M_U\) and \(j_W(j_U) = j_U\restriction M_W\).
\end{thm}

In particular \(U \I W\) and \(W\I U\).

Here we will provide a new, more general, and somewhat easier proof of this fact.

\begin{lma}\label{CommutingLemma}
Suppose \(j_0 : V\to M_0\) and \(j_1 : V\to M_1\) are ultrapower embeddings. Assume \(j_0(j_1) = j_1\restriction M_0\). Then \(j_1(j_0) = j_0\restriction M_1\).
\begin{proof}
Note that \begin{equation}\label{Models}j_1(j_0)(M_1) = j_1(M_0) = j_0(j_1)(M_0) = j_0(M_1)\end{equation} Therefore \(j_1(j_0)\) and \(j_0\restriction M_1\) are elementary embeddings from \(M_1\) into a common target model. By \cref{MinDefEmb} below, it follows that \(j_1(j_0)(\alpha) \leq j_0(\alpha)\) for all ordinals \(\alpha\). 

Let \(\xi\) be the least ordinal such that \(M_1 = H^{M_1}(j_1[V]\cup \{\xi\})\). We claim \(j_1(j_0)(\xi) = j_0(\xi)\). By the previous paragraph, it suffices to show that \(j_0(\xi) \leq j_1(j_0)(\xi)\). By the elementarity of \(j_0\), \(j_0(\xi)\) is the least ordinal \(\xi'\) such that \[j_0(M_1)= H^{j_0(M_1)}(j_0(j_1)[M_0]\cup \{\xi'\})\] To show \(j_0(\xi) \leq j_1(j_0)(\xi)\), it therefore suffices to show that \[j_0(M_1) = H^{j_0(M_1)}(j_0(j_1)[M_0]\cup \{j_1(j_0)(\xi)\})\] Note that 
\begin{align}
j_0(M_1) &= j_1(M_0)\label{Models2}\\
 &= H^{j_1(M_0)}(j_1[M_0]\cup j_1(j_0)[M_1])\label{MinCor}\\
  &= H^{j_1(M_0)}(j_1[M_0]\cup \{j_1(j_0)(\xi)\})\label{xithing}\\
  &= H^{j_0(M_1)}(j_0(j_1)[M_0]\cup \{j_1(j_0)(\xi)\})\label{AllThing}
 \end{align}
 For \cref{Models2}, we use \cref{Models}. For \cref{MinCor}, we use \cref{InternalMinimal}. For \cref{xithing}, we use \(M_1 = H^{M_1}(j_1[V]\cup \{\xi\})\). For \cref{AllThing}, we use \cref{Models} again and the fact that \(j_0(j_1) = j_1\restriction M_0\).
 
Since \[j_0\circ j_1 = j_0(j_1) \circ j_0 = j_1\circ j_0 = j_1(j_0)\circ j_1\] we have \(j_0\restriction j_1[V] = j_1(j_0) \restriction j_1[V]\). Therefore we have \[j_0 \restriction (j_1[V]\cup \{\xi\}) = j_1(j_0) \restriction (j_1[V]\cup \{\xi\})\] Since \(M_1 = H^{M_1}(j_1[V]\cup \{\xi\})\), it follows that \(j_0\restriction M_1 = j_1(j_0)\), as desired.
\end{proof}
\end{lma}

\begin{qst}
Can this lemma be proved if \(j_0\) and \(j_1\) are only assumed to be extender embeddings?
\end{qst}

Using the lemma we prove \cref{CommutingUltrapowers}.
\begin{proof}[Proof of \cref{CommutingUltrapowers}]
Let \(j_0 = j_W\) and \(j_1 = j_U\), let \(M_0 = M_W\), and let \(M_1 = M_U\). Then easily \(j_W(j_U) =  j_{j_W(U)}^{M_W} = j^{M_W}_{U} = j_U\restriction M_W\). Thus \(j_0(j_1)= j_1\restriction M_0\), so by \cref{CommutativeLemma}, \(j_1(j_0) = j_0\restriction M_1\). In other words, \(j_U(j_W) = j_W\restriction M_U\). This completes the proof.
\end{proof}
The hardest part of \cref{CommutingLemma} is to show that \(j_1(j_0)(\xi) = j_0(\xi)\), but this can be achieved much more easily under the assumptions of \cref{CommutingUltrapowers} (with \(j_0 = j_W\) and \(j_1 = j_U\)) since then \(j_0(\xi) = \xi  = j_1(j_0)(\xi)\). Using this observation, one easily obtains the following extension of \cref{CommutingUltrapowers}:

\begin{thm}
Suppose \(E\) is an extender and \(j : V\to M\) is an elementary embedding such that \(E\in V_{\textsc{crt}(j)}\) and \((j_E)^{M} = j_E\restriction M\). Then \(j_E(j) = j\restriction M_E\).\qed
\end{thm}

We remark on another characterization of commuting ultrapowers which seems to explain the term ``commuting":
\begin{prp}
Suppose \(U\) and \(W\) are countably complete uniform ultrafilters. Then the following are equivalent:
\begin{enumerate}[(1)]
\item \(j_U(j_W) = j_W\restriction M_U\).
\item \(j_U(W) = \s U W\).
\item For any set \(A\), \(\forall^U x\ \forall^W y\ (x,y)\in A\iff \forall^W y\ \forall^U x\ (x,y)\in A\).
\end{enumerate}
\begin{proof}
To see that (1) implies (2), assume (1) and note that by \cref{Pushforward}, \(\s U W\) is the ultrafilter derived from \(j_W\restriction M_U\) using \(j_W(j_U)([\text{id}]_W)\). But given (1) and applying \cref{CommutingLemma}, it follows that \(\s U W\) is the ultrafilter derived from \(j_U(j_W)\) using \(j_U([\text{id}]_W)\), which by the elementarity of \(j_U\) is equal to \(j_U(W)\).

To see that (2) implies (1), assume (2) and note that \(j_U(j_W) = j_{j_U(W)}^{M_U} = j_{\s U W}^{M_U} = j_W\restriction M_U\) by \cref{Pushforward}.

That (2) and (3) are equivalent, one applies Los's theorem. Suppose \(A\) is a set, and for any \(x\) let \(A_x = \{y \in \textsc{sp}(W) : (x,y)\in A\}\). Thus \([A_x]_U\) is a typical element of \(j_U(P(\textsc{sp}(W)))\).
On the one hand,
\begin{align*}
\forall^U x\ \forall^W y\ (x,y)\in A &\iff \forall^U x\ x : A_x\in W\\
&\iff [A_x]_U\in j_U(W)
\end{align*}
On the other hand,
\begin{align*}
\forall^W y\ \forall^U x\ (x,y)\in A &\iff \forall^W y\ (\forall^U x\ y\in A_x)\\
&\iff \forall^W y\ (j_U(y)\in [A_x]_U)\\ 
&\iff j_U^{-1}[A_x]_U\in W\\
&\iff [A_x]_U\in \s U W
\end{align*}
It follows that (2) and (3) are equivalent.
\end{proof}
\end{prp}

We now prove some ``converses" of \cref{CommutingUltrapowers}. It is easy to produce examples of uniform countably complete ultrafilters \(U\) and \(W\) such that \(U\I W\) and \(W\I U\) yet \(W\) is not \(\textsc{sp}(U)^+\)-complete and \(U\) is not \(\textsc{sp}(W)^+\)-complete. These examples are formed by iterating ultrapowers that do satisfy the hypotheses of \cref{CommutingUltrapowers}. These iterations always leave ``gaps"  in the spaces of the associated extenders:

\begin{defn}
An extender \(E\) is {\it gap-free} if the set of cardinals \(\lambda\) such that \(E\) has a generator in \([\sup j_E[\lambda],j_E(\lambda)]\) is an interval. \end{defn}
\begin{defn}
The {\it natural length} of an extender \(E\), denoted \(\nu(E)\), is the strict supremum of its generators.
\end{defn}

\begin{prp}
Suppose \(E_0\) and \(E_1\) are gap-free extenders such that \(j_{E_0}(j_{E_1}) = j_{E_1}\restriction M_{E_0}\) and \(j_{E_1}(j_{E_0}) = j_{E_0}\restriction M_{E_1}\). Then \(\nu(E_0) < \textsc{crt}(E_1)\) or \(\nu(E_1) < \textsc{crt}(E_0)\).
\begin{proof}
Let \(\kappa_0 = \textsc{crt}(E_0)\) and \(\kappa_1 = \textsc{crt}(E_1)\). Let \(I_0\) be the interval of regular cardinals at which \(E_0\) is discontinuous.  Let \(I_1\) be the interval of regular cardinals at which \(E_1\) is discontinuous.   Note that \(j_{E_0}(\kappa_1) = \kappa_1\) and \(j_{E_1}(\kappa_0) = \kappa_0\). Therefore \(\kappa_1\notin I_0\) and \(\kappa_0\notin I_1\). It follows that \(I_0\) and \(I_1\) are disjoint. Therefore either \(I_0\subseteq \kappa_1\) or \(I_1\subseteq \kappa_0\). By symmetry we may assume \(I_0\subseteq \kappa_1\).

We claim that \(I_0\) is bounded below \(\kappa_1\). To see this, suppose \(I_0\) is unbounded in \(\kappa_1\). Then \(\kappa_1 \in j_{E_1}(I_0)\). Therefore \(j_{E_1}(j_{E_0})\) is discontinuous at \(\kappa_1\), contradicting that \(j_{E_1}(j_{E_0}) = j_{E_0}\restriction M_{E_1}\) and \(j_{E_0}(\kappa_1) = \kappa_1\).

Finally, \(\nu(E_0)\leq \sup j_{E_0}[I_0] < \sup j_{E_0}[\kappa_1] = \kappa_1\), proving the proposition.
\end{proof}
\end{prp}

The following proposition shows that there are 2-cycles in the internal relation on extenders that do not arise from the commuting ultrapowers of \cref{CommutingUltrapowers}:

\begin{prp}\label{IterationCycle}
Suppose \(\kappa\) is a measurable limit of measurable cardinals. Then there is an extender \(E\) with natural length \(\kappa\) such that for any normal ultrafilter \(U\) on \(\kappa\), \(E\I U\) and \(U\I E\). Moreover \(j_U(j_E) \neq j_E\restriction M_U\) and \(j_E(j_U)\neq j_U\restriction M_E\).
\begin{proof}
Suppose \(E\) is an extender with natural length \(\kappa\) that has the property that for any \(\alpha < \kappa\), if \(i : M_{E\restriction \alpha} \to M_E\) is the canonical factor embedding, then \(i\) is definable from parameters over \(M_{E\restriction \alpha}\). (Such an extender can be constructed as a linear iteration of normal ultrafilters up to \(\kappa\).)

Let \(U\) be a normal ultrafilter on \(\kappa\). Then \(j_U(E)\restriction \kappa  = E\). So \(E\in M_U\). Moreover easily \(j_E^{M_U} = j_E\restriction M_U\) so \(E\I U\). We will show \(U\I E\) as well. Let \(i : M_E^{M_U}\to M_{j_U(E)}^{M_U}\) be the canonical factor embedding. By our assumption about \(E\), \(i\) is definable from parameters over \(M_E^{M_U}\). 

We claim that \(i\circ j_E(j_U) = j_U\). This implies \(U \I E\), since \(i\circ j_E(j_U)\) is definable over \(M_E\) using the fact that \(i\) is definable over \(M^{M_U}_E = M^{M_E}_{j_E(U)}\). We have \[i\circ j_E(j_U) \circ j_E = i\circ j_E\circ j_U = j_U(j_E) \circ j_U = j_U\circ j_E\] so \(i\circ j_E(j_U)\restriction j_E[V] = j_U\restriction j_E[V]\). Morover \(i\circ j_E(j_U)\restriction \kappa = \text{id}\restriction \kappa = j_U\restriction \kappa\). Since \(M_E = H^{M_E}(j_E[V]\cup \kappa)\), it follows that \(i\circ j_E(j_U) = j_U\restriction M_E\), as claimed.
\end{proof}
\end{prp}

Under the Ultrapower Axiom, commuting ultrapowers yield the only 2-cycles in the internal relation. To prove this, we need the following lemma, which is \cite{RF} Theorem 5.12.

\begin{defn}
Suppose \(j_0 : V\to M_0\) and \(j_1 : V\to M_1\). A pair of internal ultrapower embeddings \((i_0,i_1) : (M_0,M_1)\to N\) is a {\it comparison} of \((j_0,j_1)\) if \(i_0 \circ j_0 = i_1\circ j_1\).
\end{defn}

\begin{lma}[UA]
Any pair of ultrapower embeddings admits a unique minimal comparison.\qed
\end{lma}

\begin{lma}\label{InternalComparison}
If \(U\) and \(W\) are countably complete ultrafilters and \(U\I W\) then \((j_U(j_W),j_U\restriction M_W)\) is a minimal comparison of \((j_U,j_W)\).
\begin{proof}
The minimality of  \((j_U(j_W),j_U\restriction M_W)\) is just \cref{InternalMinimal}. The fact that \((j_U(j_W),j_U\restriction M_W)\) is a comparison is a consequence \cref{Pushforward} (which implies that \(j_U\restriction M_W\) is an internal ultrapower embedding of \(M_W\)) and the standard identity \(j_U(j_W)\circ j_U = j_U\circ j_W\).
\end{proof}
\end{lma}

\begin{thm}[UA]\label{CommutingConverse}
Suppose \(U\) and \(W\) are countably complete ultrafilters. The following are equivalent:
\begin{enumerate}[(1)]
\item \(U\I W\) and \(W\I U\).
\item \(j_U(j_W) = j_W\restriction M_U\) and  \(j_W(j_U) = j_U\restriction M_W\).
\end{enumerate}
\begin{proof}
Suppose \(U\I W\). Then \((j_U(j_W),j_U\restriction M_W)\) is the unique minimal comparison of \((j_U,j_W)\). If in addition \(W\I U\), \((j_W\restriction M_U,j_W(j_U))\) is also the unique minimal comparison of \((j_U,j_W)\). Therefore \((j_U(j_W),j_U\restriction M_W) = (j_W\restriction M_U,j_W(j_U))\), so (2) and (3) hold.
\end{proof}
\end{thm} 

\section{The generalized seed order}\label{SeedGen}
In this section we explore the relationship between the internal relation and the seed order of \cite{SO}. We actually define a somewhat more general order here called the generalized seed order. We start by defining the orders in which we will ultimately be interested:
\begin{defn}
A {\it pointed inner model} is a pair \((M,\alpha)\) where \(M\) is an inner model and \(\alpha\) is an ordinal. 

The {\it \(\Sigma_n\)-seed order} is the order on pointed inner models defined by \((M_0,\alpha_0) \swo (M_1,\alpha_1)\) if there exists \((i_0,i_1) : (M_0,M_1)\to N\) with \(i_1\) is \(\Sigma_n\)-definable over \(M_1\) from parameters and \(i_0(\alpha_0) < i_1(\alpha_1)\).
\end{defn}

For the basic analysis of the \(\Sigma_n\)-seed order, we proceed abstractly. For the time being, we fix a category \(\mathcal C\) and two collections \(J\) and \(I\) of morphisms of \(\mathcal C\). 

\begin{defn}
We say \(I\) is {\it wellfounded} if there is no sequence \(\langle i_n : n < \omega\rangle\) of elements of \(I\) such that \(\text{cod}(i_n) = \text{dom}(i_{n+1})\) for all \(n < \omega\). If \(u\) is an object of \(\mathcal C\), we say \(I\) is {\it wellfounded below \(u\)} if there is no sequence of morphisms \(\langle i_k : k <\omega\rangle\) in \(I\) such that \(\text{dom}(i_0) = u\) and for all \(k < \omega\), \(\text{cod}(i_k) = \text{dom}(i_{k+1})\).
\end{defn}

\begin{defn}
We say \((J, I)\) has the {\it shift property} if for any \(j : v \to w\) in \( J\) and any \(i : v\to u\) in \(I\), there is an object \(x\in \mathcal C\) admitting morphisms  \(i' : w\to x\) in \(I\) and \(j' : u \to x\) in \(J\).
\end{defn}

\begin{defn}
The {\it \((J,I)\)-seed order} is defined on the objects of \(\mathcal C\) by setting \(u < w\) if there is some \(x\in \mathcal C\) admitting morphisms \(j : u\to x\) in \(J\) and \(i : w \to x\) in \(I\).
\end{defn}

\begin{lma}\label{Transitivity}
Suppose \(J\) and \(I\) are closed under composition and \((J,I)\) has the shift property. Then the \((J,I)\)-seed order is transitive.
\begin{proof}
Fix \(u_0,u_1,u_2\in \mathcal C\) with \(u_0 < u_1 < u_2\). We must show \(u_0 < u_2\). 

Since \(u_0 < u_1\), we can find morphisms \(j_0 : u_0\to w_0\) in \(J\) and \(i_1 : u_1 \to w_0\) in \(I\). Since \(u_1 < u_2\) we can find morphisms \(j_1 : u_1\to w_1\) in \(J\) and \(i_2 : u_2 \to w_1\) in \(I\). By the shift property applied to \(i_1\) and \(j_1\), we can find \(x\in \mathcal C\) admitting \(j' : w_0\to x\) in \(J\) and \(i' : w_1\to x\) in \(I\). Then by the closure of \(J\) and \(I\) under composition, \(j'\circ j_0 : u_0 \to x\) is in \(J\) and \(i'\circ i_2 : u_2\to x\) is in \(I\), so \(u_0 < u_2\) as desired.
\end{proof}
\end{lma}

\begin{thm}\label{GeneralizedWF}
Suppose \(I\) is wellfounded and \((J,I)\) has the shift property. Then the \((J,I)\)-seed order is wellfounded.
\begin{proof}
We start with a simple construction. Given a \((J,I)\)-seed order descending sequence \(u_0 > u_1 > u_2 > \cdots\), we show how to produce another such sequence \(u^*_0 > u^*_1 > u^*_2 > \cdots\) such that there is a morphism \(i_0 : u_0\to u_0^*\) in \(I\). Since \(u_n > u_{n+1}\), we can fix an object \(u_n^*\) and morphisms \(i_n : u_n\to u^*_n\) in \(I\) and \(j_n : u_{n+1}\to u^*_n\) in \(J\). We claim that for \(n <\omega\), \(u_n^* > u_{n+1}^*\). To see this, use the shift property on \(i_{n+1} : u_{n+1}\to u^*_{n+1}\) and \(j_n : u_{n+1}\to u^*_n\) to obtain an object \(w\in \mathcal C\) admitting morphisms \(j^*_n : u^*_{n+1}\to w\) in \(J\) and \(i^*_n : u^*_n \to w\) in \(I\). The existence of these morphisms implies \(u^*_n > u^*_{n+1}\).

Now assume towards a contradiction that the \((J,I)\)-seed order is illfounded. Fix \(u_0^0 > u^0_1 > \cdots\). By recursion we define objects \(\{u^n_m :  n ,m < \omega\}\) of \(\mathcal C\) and morphisms \(i^n : u^n_0\to u^{n+1}_0\) in \(I\). Suppose \(u^n_0 > u^n_1 > \cdots\) has been defined. By the previous paragraph we can produce \(u^{n+1}_0 > u^{n+1}_1 > \cdots\) and a morphism \(i^n :u^n_0\to u^{n+1}_0\) in \(I\). The sequence \(\langle i^n : n < \omega\rangle\) contradicts the wellfoundedness of \(I\).
\end{proof}
\end{thm}

\begin{cor}\label{Strictness}
Suppose \(I\) is wellfounded and \((J,I)\) has the shift property. Then for any objects \(u_0,u_1\in \mathcal C\), one of \(\hom(u_0,u_1)\cap I\) and \(\hom(u_0,u_1)\cap J\) is empty.
\begin{proof}
Otherwise \(u_0 < u_0\) in the \((J,I)\)-seed order, contradicting \cref{GeneralizedWF}.
\end{proof}
\end{cor}

We now apply these general facts to a specific category:

\begin{defn}
We denote by \(\mathscr C\) the category of pointed inner models with all elementary embeddings.
\end{defn}

\(\mathscr C\) is a pretty large category, but everything we do can be formalized quite easily in NBG. The following lemma is really a schema that is proved for each fixed natural number \(n\) in the metatheory.

\begin{lma}\label{Example}
Let \(J\) be the collection of morphisms \(j : (M,\alpha)\to (N,\beta)\) in \(\mathscr C\) such that \(j(\alpha) = \beta\). Let \(I\) be the collection of morphisms \(i : (M,\alpha)\to (N,\beta)\) in \(\mathscr C\) such that \(i\) is \(\Sigma_n\)-definable from parameters over \(M\) and \(i(\alpha) > \beta\). Then \((J,I)\) has the shift property and \(I\) is wellfounded.
\begin{proof}
We first prove the shift property. Suppose \(j : (M,\alpha)\to (N,j(\alpha))\) is in \(J\) and \(i : (M,\alpha)\to (P,\beta)\) is in \(I\). Let \(Q = j(P)\) and \(\gamma = j(\beta)\). Let \(i' = j(i)\) and let \(j' = j\restriction P\). Obviously \(j' : (P,\beta)\to (Q,\gamma)\) is in \(C\). To finish, we just need to show that \(i' : (N,j(\alpha))\to (Q,\gamma)\) belongs to \(I\). First, \(j(i)\) is \(\Sigma_n\)-definable over \(N\) using the definition of \(i\) with its parameters shifted by \(j\). Second \(i'(j(\alpha)) = j(i(\alpha)) > j(\beta) = \gamma\). Thus  \(i' : (N,j(\alpha))\to (Q,\gamma)\) is in \(I\).

We finally show that \(I\) is wellfounded. This follows from Kunen's proof \cite{KunenLU} of the wellfoundedness of iterated ultrapowers. We sketch this argument here. We require the following claim, which is proved by an easy absoluteness argument.

\begin{clm}\label{Absoluteness}
Suppose \((M,\alpha)\) is a pointed inner model and \(I\) is illfounded below \((M,\alpha)\). Then \((M,\alpha)\) satisfies that \(I\) is illfounded below \((M,\alpha)\).\qed
\end{clm}

Suppose towards a contradiction that \(I\) is illfounded. Fix an inner model \(M\) such that for some ordinal \(\alpha\), \(I\) is illfounded below \((M,\alpha)\). Let \(\alpha\) be the least ordinal such that \(I\) is illfounded below \((M,\alpha)\). By \cref{Absoluteness}, \(M\) satisfies that \(\alpha\) is the least ordinal \(\alpha'\) such that \(I\) is illfounded below \((M,\alpha')\). 

Fix a sequence \(\langle i_k : k < \omega\rangle\) of elements of \(I\) with \(\text{dom}(i_0) = (M,\alpha)\) and \(\text{cod}(i_k) = \text{dom}(i_{k+1})\) for all \(k < \omega\). Let \((N,\beta) = \text{cod}(i_0)\). By the elementarity of \(i_0\), \(N\) satisfies that \(i_0(\alpha)\) is the least ordinal \(\alpha'\) such that \(I\) is illfounded below \((N,\alpha')\). By \cref{Absoluteness}, \(i_0(\alpha)\) actually is the least ordinal \(\alpha'\) such that \(I\) is illfounded below \((N,\alpha')\). But \(\beta < i_0(\alpha)\) by the definition of \(I\), and \(\langle i_k : 1\leq k < \omega\rangle\) witnesses that \(I\) is illfounded below \((N,\beta)\). This is a contradiction.
\end{proof}
\end{lma}

\begin{cor}\label{WF}
The \(\Sigma_n\)-seed order is a wellfounded strict partial order of the collection of pointed inner models.
\begin{proof}
The \(\Sigma_n\)-seed order is the \((J,I)\)-seed order on \(\mathscr C\) where \(J\) and \(I\) are as in \cref{Example}. Therefore the corollary follows from \cref{Transitivity} and \cref{GeneralizedWF}.
\end{proof}
\end{cor}

As a consequence, we have the following theorem which is often useful:

\begin{thm}\label{MinDefEmb}
Suppose \(M\) is an inner model and \(i,j : M\to N\) are elementary embeddings. Assume \(i\) is definable over \(M\) from parameters. Then for all ordinals \(\alpha\), \(i(\alpha)\leq j(\alpha)\).
\begin{proof}
Suppose there is a counterexample such that \(i\) is a \(\Sigma_n\)-definable elementary embedding. Then \(i,j : (M,\alpha)\to (M,j(\alpha))\) are morphisms of \(\mathscr C\) with \(i\in I\) and \(j\in J\), where \(I\) and \(J\) are defined as in the statement of \cref{Example}. This contradicts \cref{Strictness}.
\end{proof}
\end{thm}

\subsection{Steel's conjecture}
In this section we put down some corollaries of \cref{GeneralizedWF} for the internal relation.
\begin{thm}\label{Descent}
For any ordinal \(\alpha\), the internal relation is wellfounded on extenders discontinuous at \(\alpha\).
\end{thm}

The proof uses the following lemma:
\begin{lma}\label{DescentLemma}
Assume \(E\I F\) are extenders and \(\alpha\) is an ordinal at which \(j_E\) is discontinuous. Then \[(M_E,\sup j_E[\alpha]) \swo (M_F,\sup j_F[\alpha])\] in the \(\Sigma_2\)-seed order. 
\begin{proof}
As usual, consider the comparison \((j_E(j_F),j_E) : (M_E,M_F)\to j_E(M_F)\):
\begin{align*}
 j_E(j_F)(\sup j_E[\alpha]) &< \sup j_E(j_F)[j_E(\alpha)]\\
 &= j_E(\sup j_F[\alpha]) 
\qedhere
\end{align*}
\end{proof}
\end{lma}

\begin{proof}[Proof of \cref{Descent}]
Assume \(E_0 \gI E_1 \gI E_2\gI\cdots\) is a descending sequence of extenders in the internal relation such that for all \(n < \omega\), \(j_{E_n}\) is discontinuous at \(\alpha\). By \cref{DescentLemma}, \[(M_{E_0},\sup j_{E_0}[\alpha]) \slwo (M_{E_1},\sup j_{E_1}[\alpha])\slwo (M_{E_2},\sup j_{E_2}[\alpha])\slwo \cdots\] contrary to the wellfoundedness of the \(\Sigma_2\)-seed order, \cref{WF}.
\end{proof}

As a corollary, we have some additional information about 2-cycles in the internal relation:
\begin{cor}
Suppose \(E\) and \(F\) are extenders such that \(E\I F\) and \(F\I E\). Then \(E\) and \(F\) have no common discontinuity points.
\begin{proof}
If \(E\) and \(F\) are both discontinuous at \(\alpha\), then \(E\gI F\gI E\gI F\gI \cdots\) witnesses the illfoundedness of the internal relation on extenders discontinuous at \(\alpha\), contradicting \cref{Descent}.
\end{proof}
\end{cor}

We now use \cref{Descent} to prove a conjecture of Steel \cite{SteelWellfounded}.
\begin{thm}\label{SteelConj}
Suppose \(E_0\gmo E_1\gmo E_2\gmo \cdots\) is a sequence of rank-to-rank extenders of length \(\lambda\). Then \(\{ \textsc{crt}(E_n) :  n < \omega\}\) is cofinal in \(\lambda\).
\end{thm}

We need the following lemma which appears as part of the proof of \cite{SteelWellfounded} Theorem 2.2:
\begin{lma}
Suppose \(E\) and \(F\) are rank-to-rank extenders of length \(\lambda\). If \(E\mo F\) then \(j_E^{M_F} = j_E\restriction M_F\), so \(E\I F\).\qed
\end{lma}

\begin{proof}[Proof of \cref{SteelConj}]
Let \(\bar \lambda = \sup_{n <\omega} \textsc{crt}(E_n)\). Suppose towards  a contradiction that \(\bar \lambda < \lambda\). Note that for all \(n < \omega\), \(j_{E_n}\) is discontinuous at every regular cardinal in \([\textsc{crt}(E_n),\lambda]\). Therefore fix any regular cardinal \(\delta\) in the interval \([\bar \lambda,\lambda]\). Then for all \(n\), \(j_{E_n}\) is discontinuous at \(\delta\). The sequence \(E_0\gI E_1\gI E_2\gI\cdots\) therefore contradicts the wellfoundedness of the internal relation of extenders discontinuous at \(\delta\).
\end{proof}
\section{UA and the internal relation}
We now take a closer look at the structure of the internal relation on countably complete ultrafilters assuming UA.

\subsection{Translation functions and the internal relation}
\begin{defn}
Suppose \(U\) is a countably complete ultrafilter and in \(M_U\), \(W'\) is a countably complete uniform ultrafilter on an ordinal \(\delta'\). Then the {\it \(U\)-limit of \(W'\)} is the ultrafilter \[U^-(W') = \{X\subseteq \delta : j_U(X)\cap \delta'\in W'\}\] where \(\delta\) is the least ordinal such that \(\delta'\leq j_U(\delta)\).
\end{defn}

\begin{defn}[UA]
For \(U,W\in \Un\), \(\tr{W}{U}\) denotes the \(\swo^{M_W}\)-least \(U'\in \Un^{M_W}\) such that \(W^-(U') = U\). The function \(t_W : \Un\to \Un^{M_W}\) is called the {\it translation function} associated to \(W\).
\end{defn}

An immediate consequence of the definition of translation functions is the following bound:

\begin{lma}[UA]\label{BoundingLemma}
For any \(U,W\in \Un\), \(\tr U W \wo^{M_U} j_U(W)\).
\begin{proof}
Note that \(U^-(j_U(W)) = W\), so by the minimality of \(\tr U W\), \(\tr U W \wo^{M_U} j_U(W)\).
\end{proof}
\end{lma}

We will use the following theorem, which appears as \cite{SO} Theorem 5.3.

\begin{thm}[UA]\label{Reciprocity1}
Suppose \(U_0,U_1\in \Un\) and \((i_0,i_1) : (M_{U_0},M_{U_1})\to N\) is a comparison of \((j_{U_0},j_{U_1})\). Then \(\tr {U_0} {U_1}\) is the ultrafilter derived from \(i_0\) using \(i_1([\textnormal{id}]_{U_1})\).
\end{thm}

As a corollary we obtain the following information about the relationship between translation functions and the internal relation, generalizing \cref{CommutingConverse}.

\begin{thm}[UA]\label{InternalTranslation}
Suppose \(U\) and \(W\) are uniform countably complete ultrafilters. Then the following are equivalent:
\begin{enumerate}[(1)]
\item \(U\I W\).
\item \(\tr W U = \s W U\).
\item \(\tr U W = j_U(W)\).
\item \(j_U(W)\wo^{M_U} \tr U W\).
\end{enumerate}
\begin{proof}
To see (1) implies (2) and (3), assume (1). By \cref{InternalComparison}, \((j_U(j_W),j_U\restriction M_W)\) is a comparison of \((j_U,j_W)\). By \cref{Pushforward}, \(\s W U\) is the \(M_W\)-ultrafilter derived from \(j_U\restriction M_W\) using \(j_U(j_W)([\text{id}]_U)\) so (2) holds. By the elementarity of \(j_U\), \(j_U(W)\) is the ultrafilter derived from \(j_U(j_W)\) using \(j_U([\text{id}]_W)\), so (3) holds.

(2) implies (1) by \cref{Pushforward}.

We now show (3) implies (1). Let \((i_0,i_1) : (M_U,M_W)\to N\) be the unique minimal comparison of \((j_U,j_W)\). We claim \(i_1 = j_U\restriction M_W\). By the minimality of \((i_0,i_1)\), \(N = H^N(i_0[M_U]\cup i_1[M_W]) = H^N(i_1[M_U]\cup \{i_1([\text{id}]_W)\})\). Therefore \(i_0\) is the ultrapower by the ultrafilter derived from \(i_0\) using \(i_1([\text{id}]_W)\). But the ultrafilter derived from \(i_0\) using \(i_1([\text{id}]_W)\) is \(j_U(W)\) by \cref{Reciprocity1} and the assumption that (3) holds. Therefore  \(i_0 = j_U(j_W)\) and \(i_1([\text{id}]_W) = j_U([\text{id}]_W)\). Moreover \(i_1 \circ j_W = i_0 \circ j_U = j_U(j_W) \circ j_U = j_U \circ j_W\). Therefore \[i_1 \restriction j_W[V]\cup \{[\text{id}]_W\} = j_U \restriction j_W[V]\cup \{[\text{id}]_W\}\] Since \(M_W = H^{M_W}(j_W[V]\cup \{[\text{id}]_W\})\), it follows that \(i_1 = j_U\restriction M_W\), as desired.

Finally (3) implies (4) trivially and (4) implies (3) by \cref{BoundingLemma} and the antisymmetry of the seed order.
\end{proof}
\end{thm}

\subsection{\(j\) on the ordinals}
In this subsection we study the relationship between the internal relation and the action of ultrapower embeddings on the ordinals.

The first thing we show is that whether \(U\I W\) really only depends on the fixed points of \(j_U\).

\begin{prp}[UA]\label{FixThm}
Suppose \(U_0\) and \(U_1\) are countably complete ultrafilters such that \(j_{U_0}\) fixes every ordinal fixed by \(j_{U_1}\). Then for any countably complete ultrafilter \(W\) with \(U_1\I W\), \(U_0\I W\).
\end{prp}

To prove this we use an analysis of the seed order on pointed ultrapowers. 
\begin{defn}
A {\it pointed ultrapower} \(\mathcal M\) is a pair \(\mathcal M = (M,\alpha)\) where \(M\) is an ultrapower of \(V\) and \(\alpha\) is an ordinal. The collection of pointed ultrapowers is denoted by \(\mathscr P\). If If \(\mathcal M\) is a pointed ultrapower, then \(\alpha_\mathcal M\) denotes the ordinal \(\alpha\) such that \(\mathcal M = (M,\alpha)\) for some inner model \(M\).
\end{defn}

If \(\mathcal M = (M,\alpha)\), we will abuse notation by writing \(\mathcal M\) when we really mean the inner model \(M\). 

\begin{defn}
The {\it completed seed order} is defined on \(\mathcal M_0,\mathcal M_1\in \mathscr P\) by setting \(\mathcal M_0 \swo \mathcal M_1\) if there exists \((i_0,i_1) : (\mathcal M_0,\mathcal M_1)\to N\) such that \(i_1\) is an internal ultrapower embedding of \(\mathcal M_1\) and \(i_0(\alpha_{\mathcal M_0}) < i_1(\alpha_{\mathcal M_1})\).
\end{defn}

For any \(n\geq 2\), the completed seed order is the restriction of the \(\Sigma_n\)-seed order to the collection of pointed ultrapowers. Therefore it is a wellfounded strict partial order. It is not true, however, that UA implies that the completed seed order is total. To explain this a bit more clearly, it is worth introducing the following nonstrict version of the completed seed order:

\begin{defn}
The {\it nonstrict completed seed order} is defined on \(\mathcal M_0,\mathcal M_1\in \mathscr P\) by setting \(\mathcal M_0 \wo \mathcal M_1\) there exists \((i_0,i_1) : (\mathcal M_0,\mathcal M_1)\to N\) such that \(i_1\) is an internal ultrapower embedding of \(\mathcal M_1\) and \(i_0(\alpha_{\mathcal M_0}) \leq i_1(\alpha_{\mathcal M_1})\). We write \(\mathcal M_0\equiv_S \mathcal M_1\) if \(\mathcal M_0\wo \mathcal M_1\) and \(\mathcal M_1\wo \mathcal M_0\). 
\end{defn}

If \(\mathcal M_0\equiv_S\mathcal M_1\), then \(\mathcal M_0\) and \(\mathcal M_1\) are incomparable in \(\swo\). Many instances of \(\equiv_S\) on pointed ultrapowers arise from the following trivial lemma:

\begin{lma}\label{EquivLemma}
If \(\mathcal M\) and \(\mathcal N\) are pointed ultrapowers and \(j : \mathcal M\to \mathcal N\) is an internal ultrapower embedding with \(j(\alpha_\mathcal M) = \alpha_\mathcal N\), then \(\mathcal M\equiv_S\mathcal N\).\qed
\end{lma}

The question of whether all instances of \(\equiv_S\) arise from the previous lemma remains open:

\begin{qst}[UA]
Suppose \(\mathcal M_0,\mathcal M_1\in \mathscr P\) and \(\mathcal M_0 \equiv \mathcal M_1\). Does there exist \(\mathcal N\in \mathscr P\) admitting internal ultrapowers \(j_0 : \mathcal N\to \mathcal M_0\) and \(j_1 : \mathcal N\to \mathcal M_1\) with \(j_0(\alpha_\mathcal N) = \alpha_{\mathcal M_0}\) and \(j_1(\alpha_\mathcal N) = \alpha_{\mathcal M_1}\)?
\end{qst}

This is related to the question of characterizing greatest lower bounds in the Rudin-Frolik order, which are proved to exist in \cite{RF} Theorem 7.3 without being described explicitly. In particular, an affirmative answer to this question is equivalent to the statement that for any countably complete ultrafilters \(U_0\) and \(U_1\), the ultrapower by the greatest lower bound of \(U_0\) and \(U_1\) is isomorphic to \(i_0[M_{U_0}]\cap i_1[M_{U_1}]\) for any comparison \((i_0,i_1)\) of \((j_{U_0},j_{U_1})\). (In particular, this would imply the distributivity of the Rudin-Frolik lattice under UA.)

In any case, the nonstrict completed seed order is a prewellorder under UA:

\begin{lma}[UA]
The nonstrict completed seed order prewellorders \(\mathscr P\). In fact, for any \(\mathcal M_0, \mathcal M_1\in \mathscr P\), either \(\mathcal M_0\swo \mathcal M_1\) or \(\mathcal M_1\wo \mathcal M_0\). Moreover for \(\mathcal M_0,\mathcal M_1\in \mathscr P\), \(\mathcal M_0\wo \mathcal M_1\) if and only if there are internal ultrapowers \((i_0,i_1) : (\mathcal M_0,\mathcal M_1)\to N\) with \(i_0(\alpha_{\mathcal M_0})\leq i_1(\alpha_{\mathcal M_1})\).\qed
\end{lma}

The completed seed order completes the seed order in the following sense.
\begin{defn}
Define \(\Phi :\Un \to \mathscr P\) by \(\Phi(U) = (M_U,[\textnormal{id}]_U)\)
\end{defn}

Assuming just ZFC, it is not clear that \(\Phi(U)\swo \Phi(W)\) implies \(U\sE W\), but this is a consequence of UA (or \(V = \text{HOD}\)).

\begin{lma}[UA]\label{Completion}
\(\Phi\) is an order embedding from \((\Un,\swo)\) into \((\mathscr P,\swo)\).\qed
\end{lma}

The rank of an ultrafilter in the completed seed order may not exist since the completed seed order may not be setlike. The following theorem (\cite{Frechet} Theorem 11.16) characterizes this behavior:

\begin{thm}[UA]
Exactly one of the following holds:
\begin{enumerate}[(1)]
\item The completed seed order is setlike.
\item There is a supercompact cardinal.\qed
\end{enumerate}
\end{thm}

We therefore consider restricted versions of the completed seed order in order to obtain rank functions that take values in the ordinals.

\begin{defn}
For any cardinal \(\delta\), \(\mathscr P_{\leq\delta}\) denotes the class of pointed ultrapowers \((M,\alpha)\) such that \(M = M_U\) for some \(U\in \Un_{\leq\delta}\).

For any pointed ultrapower \(\mathcal M\), we denote the rank of \(\mathcal M\) in the completed seed order on \(\mathscr P_{\leq\delta}\) by \(|\mathcal M|_{\leq\delta}\). For \(U\in \Un\), we let \(|U|_{\leq \delta} = |\Phi(U)|_{\leq\delta}\).
\end{defn}

If \(\mathcal M\) happens to have an ordinal rank in the completed seed order on \(\mathscr P\), then this rank is equal to the eventual value of \(|\mathcal M|_{\leq\delta}\) for \(\delta\) arbitrarily large.

Towards \cref{FixThm}, we show the following fact:

\begin{thm}[UA]\label{FixThm2}
Suppose \(U,W\in \Un_{\leq\delta}\). Then \(U\I W\) if and only if \(j_U\) fixes \(|W|_{\leq\delta}\).
\end{thm}

This requires some lemmas. The first is an absoluteness property which sets the completed seed order apart from the usual seed order:

\begin{lma}[UA]\label{PointedAbs}
Suppose \(\delta\) is a cardinal and \(U\in \Un_{\leq\delta}\). Then \(\mathscr P^{M_U}_{\leq j_U(\delta)}\subseteq\mathscr P_{\leq \delta}\). Moreover for any \(\mathcal M\in \mathscr P_{\leq \delta}\), there is some \(\mathcal M' \in\mathscr P^{M_U}_{\leq j_U(\delta)}\) with \(\mathcal M \equiv_S \mathcal M '\). Therefore for any \(\mathcal M \in \mathscr P_{\leq j_U(\delta)}^{M_U}\), \(|\mathcal M|_{\leq \delta} = |\mathcal M|_{\leq j_U(\delta)}^{M_U}\).
\begin{proof}
The fact that \(\mathscr P^{M_U}_{\leq j_U(\delta)}\subseteq\mathscr P_{\leq \delta}\) amounts to the standard fact that an iterated ultrapower that hits ultrafilters on ordinals less than or equal to the image of \(\delta\) is given by a single ultrapower by an ultrafilter on \(\delta\).

To see that for any \(\mathcal M\in \mathscr P_{\leq \delta}\), there is some \(\mathcal M' \in\mathscr P^{M_U}_{\leq j_U(\delta)}\) with \(\mathcal M \equiv_S \mathcal M '\), fix an ultrafilter \(W\in \Un_{\leq \delta}\) such that \(\mathcal M  = M_W\). By \cref{BoundingLemma}, \(\tr U W\in \Un_{\leq j_U(\delta)}^{M_U}\). Moreover there is an internal ultrapower \(j : M_W\to M^{M_U}_{\tr U W}\). Letting \(\mathcal M ' = (M^{M_U}_{\tr U W},j(\alpha_\mathcal M))\), this implies \(\mathcal M'\in \mathscr P_{\leq \delta}\) and \(\mathcal M \equiv_S \mathcal M'\) by \cref{EquivLemma}.
\end{proof}
\end{lma}

\begin{lma}[UA]\label{TranslationRank}
Suppose \(U,W\in \Un_{\leq\delta}\). Then \(|\tr U W|^{M_U}_{\leq j_U(\delta)} = |W|_{\leq\delta}\).
\begin{proof}
Note that there is an internal ultrapower embedding \(i : M_W\to M_{\tr U W}^{M_U}\) such that \(i([\text{id}]_W) = [\text{id}]_{\tr U W}^{M_U}\). Therefore \[|W|_{\leq\delta} = |(M_{\tr U W}^{M_U},[\text{id}]_{\tr U W})|_{\leq\delta} = |(M_{\tr U W}^{M_U},[\text{id}]_{\tr U W})|_{\leq j_U(\delta)}^{M_U} = |\tr U W|^{M_U}_{\leq j_U(\delta)} \qedhere\]
\end{proof}
\end{lma}

We can finally prove \cref{FixThm2}:

\begin{proof}[Proof of \cref{FixThm2}]
Suppose first that \(U\I W\). Then \[\Phi(W) = (M_W,[\text{id}]_W)\equiv_S (j_U(M_W),j_U([\text{id}]_W)) = \Phi^{M_U}(j_U(W))\] and hence \(|W|_{\leq \delta} = j_U(|W|_{\leq\delta})\).

Conversely assume \(|W|_{\leq \delta} = j_U(|W|_{\leq\delta})\). Then \[|j_U(W)|^{M_U}_{\leq j_U(\delta)} = j_U(|W|_{\leq\delta}) = |W|_{\leq \delta} = |\tr U W|_{\leq j_U(\delta)}^{M_U}\] Therefore in \(M_U\), \(\Phi(j_U(W)) \equiv_S \Phi(\tr U W)\). By \cref{Completion}, \(j_U(W) = \tr U W\). Therefore by \cref{InternalTranslation}, \(U \I W\).
\end{proof}

As a corollary, we prove \cref{FixThm}:

\begin{proof}[Proof of \cref{FixThm}]
We may assume without loss of generality that \(U\) and \(W\) are uniform ultrafilters. Take \(\delta\) large enough that \(U_0,U_1,\) and \(W\) belong to \(\Un_{\leq\delta}\). Then by \cref{FixThm2}, since \(U_1\I W\), \(j_{U_1}(|W|_{\leq\delta}) = |W|_{\leq \delta}\). Hence by assumption \(j_{U_0}(|W|_{\leq\delta}) = |W|_{\leq \delta}\). Applying the other direction of \cref{FixThm2}, \(U_0\I W\).
\end{proof}

We now prove a theorem that gives rise to some fairly interesting examples of countably complete ultrafilters whose ultrapowers have the same fixed points. For example, we will show that there are distinct ultrapower embeddings with the same action on the ordinals.

\begin{thm}[UA]\label{rangefix}
Suppose \(U\in \textnormal{Un}\) and \(W\in \textnormal{Un}^{M_U}\). Suppose \[j_U(W)\I^{M_{U^2}} j_U(j_U)(W)\] Then \(j^{M_U}_W\circ j_U\restriction \textnormal{Ord} = j_U\restriction \textnormal{Ord}\).
\end{thm}

For the proof we need the following fact, which generalizes a lemma of Kunen. (See \cite{Larson} Lemma 1.1.26, but note that the hypothesis that \(\mu\) is a \(\kappa\)-complete ultrafilter {\it on \(\kappa\)} should have been included in the statement of that lemma).

\begin{lma}[UA]\label{KunenClaim}
Suppose \(\langle W_n : n < \omega\rangle\) is a sequence of ultrafilters such that for all \(m \geq 0\), for all \(n > m\), \(W_n \I W_m\). Then for any ordinal \(\alpha\), for all sufficiently large \(n\), \(j_{W_n}(\alpha) = \alpha\).
\begin{proof}
Suppose not, and let \(\alpha\) be the least ordinal at which \cref{KunenClaim} fails. Fix a sequence \(\langle W_n : n < \omega\rangle\) such that for all \(m\), for all \(n > m\), \(W_n \I W_m\) yet for infinitely many \(n< \omega\), \(j_{W_n}(\alpha) > \alpha\). By passing to a subsequence we may of course assume that for all \(n < \omega\), \(j_{W_n}(\alpha) > \alpha\). 

By elementarity, in \(M_{W_0}\), \(j_{W_0}(\alpha)\) is the least ordinal at which \cref{KunenClaim} fails. In particular since \(\alpha < j_{W_0}(\alpha)\), \cref{KunenClaim} holds at \(\alpha\) inside \(M_{W_0}\).

\begin{clm} For \(n \geq 1\), let \(W'_n = s_{W_0}(W_n)\). Then in \(M_{W_0}\), for all \(m\geq 1\), for all \(n > m\), \(W'_n \I W'_m\).
\begin{proof}
Let \(\lambda = \sup \{\textsc{sp}(W_n) : n < \omega\}\). Let \(\xi_m = |W_m|_\lambda\). By \cref{InternalTranslation}, \(W'_m = \tr {W_0} {W_m}\). Therefore by \cref{TranslationRank}, \(|W'_m|_{j_{W_0}(\lambda)}^{M_{W_0}} = |W_m|_{\lambda} =\xi_m\). If \(n > m\geq 1\), then \(j_{W_n}(\xi_m) =\xi_m\) by \cref{FixThm2}. But since \((j_{W_n'})^{M_{W_0}} = j_{W_n}\restriction M_{W_0}\) by \cref{Pushforward}, we have \(j_{W_n'}^{M_{W_0}} (\xi_m) = j_{W_n}(\xi_m) = \xi_m\). Since \(W_n'\in \Un_{\leq j_{W_0}(\lambda)}^{M_{W_0}}\) and since \(\xi_m = |W'_m|_{j_{W_0}(\lambda)}^{M_{W_0}}\), \cref{FixThm2} implies that \(W_n'\I W_m'\) in \(M_{W_0}\).
\end{proof}
\end{clm}

Since \(j^{M_{W_0}}_{W'_n}(\alpha) = j_{W_n}(\alpha) > \alpha\) for all \(n \geq 1\),  \(\langle W_n' : 1\leq n < \omega\rangle\) witnesses that \cref{KunenClaim} fails at \(\alpha\) in \(M_{W_0}\). This is a contradiction.
\end{proof}
\end{lma}

\begin{proof}[Proof of \cref{rangefix}]
Let \(\delta = \textsc{sp}(U)\). Let \(\langle W_\alpha : \alpha < \delta\rangle\) represent \(W\) in \(M_U\). The statement that \(j_U(W)\I^{M_{U^2}} j_U(j_U)(W)\) is equivalent to the statement that for \(U\)-almost every \(\alpha < \delta\), for \(U\)-almost every \(\beta < \delta\), \(W_\beta \I W_\alpha\). 

Fix \(\xi\in \text{Ord}\) and let us show \(j^{M_U}_W(j_U(\xi)) = j_U(\xi)\).
Suppose towards a contradiction that \(j^{M_U}_W(j_U(\xi)) > j_U(\xi)\). Thus by Los's theorem, for \(U\)-almost all \(\alpha < \delta\), \(j_{W_\alpha}(\xi) > \xi\). We denote the set of such \(\alpha\) by \(X\subseteq \delta\). 

We now construct a sequence of ordinals \(\langle \alpha_n : n < \omega\rangle\) by induction such that for all \(n\), \(j_{W_{\alpha_n}}(\xi) > \xi\) and for all \(m < n\), \(W_{\alpha_n}\I W_{\alpha_m}\). For \(\alpha < \delta\), let \[A_\alpha= \{\beta < \delta : W_\beta \I W_\alpha\}\] Then for \(U\)-almost all \(\alpha\), \(A_\alpha\in U\). Suppose \(\alpha_m\) has been defined for \(m < n\) in such a way that \(A_{\alpha_m}\in U\). We then choose \(\alpha_n \in \bigcap_{m < n} A_{\alpha_m}\) such that \(A_{\alpha_n}\in U\) and \(\alpha_n\in X\). Such an ordinal exists since \(U\)-almost all \(\alpha <\delta\) satisfy these requirements. This ensures that \(j_{W_{\alpha_n}}(\xi) > \xi\) and for all \(m < n\), \(W_{\alpha_n}\I W_{\alpha_m}\). Moreover since \(A_{\alpha_n}\in U\), we can continue the recursion.

The existence of the sequence \(\langle W_{\alpha_n} : n < \omega\rangle\) is prohibited by \cref{KunenClaim}, so we have a contradiction.
\end{proof}

The use of UA is minimal here, and there is a ZFC fact that covers the interesting cases and more:

\begin{defn}
The extenders \(E_0\) and \(E_1\) {\it commute} if \(j_{E_0}(j_{E_1}) = j_{E_1}\restriction M_{E_0}\) and \(j_{E_1}(j_{E_0}) = j_{E_0}\restriction M_{E_1}\).
\end{defn}

\begin{thm}\label{rangefix2}
Suppose \(F\) is an extender, \(E\in M_F\) is an \(M_F\)-extender, and \(j_F(j_F)(E)\) and \(j_F(E)\) commute in \(M_{F^2}\). Then \[j^{M_F}_E\circ j_F\restriction\textnormal{Ord} = j_F\restriction \textnormal{Ord}\]
\begin{proof}
We first reduce to the case that \(F\) is an ultrafilter. 

Let \(\bar F\) be the ultrafilter derived from \(F\) using \(E\). Let \[k : M_{\bar F} \to M_F\] be the factor embedding with \(k\circ j_{\bar F} = j_F\). Let \(\bar E\) be such that \(k(\bar E) = E\). It suffices to show that \(j^{M_{\bar F}}_{\bar E}\) fixes every ordinal in the range of \(j_{\bar F}\), since then for any ordinal \(\xi\), \[j^{M_F}_E(j_F(\xi)) = k(j^{M_{\bar F}}_{\bar E}(j_{\bar F}(\xi))) = k(j_{\bar F}(\xi)) = j_F(\xi)\]
\begin{clm} \(j_{\bar F}(j_{\bar F})(\bar E)\) and \(j_{\bar F}(\bar E)\) commute in \(M_{\bar F^2}\).\end{clm}
\begin{proof} Let \(k_* = j_{\bar F}(k)\), so \[k_* : M_{\bar F^2} \to M^{M_{\bar F}}_{j_{\bar F}(\bar E)}\] Let \(i = k\circ k_*\). Then \(i : M_{\bar F^2} \to M_{F^2}\). We show that \(i(j_{\bar F}(\bar E)) = j_F(E)\) and \(i(j_{\bar F}(j_{\bar F})(\bar E)) = j_F(j_F)(E)\). Since \(j_F(E)\) and \(j_F(j_F)(E)\) commute in \(M_{F^2}\), the claim then follows from the elementarity of \(i\).

This is a routine diagram chase which is easier done than said. We recommend drawing the embeddings and checking it yourself.
First,
\begin{align*}
i(j_{\bar F}(\bar E)) &= (k\circ k_*)(j_{\bar F}(\bar E))\\
&= k(j_{\bar F}(k)(j_{\bar F}(\bar E))) \\
&= k(j_{\bar F}(k(\bar E)))\\
&= k(j_{\bar F}(E))\\
&= j_F(E)
\end{align*}
Second,
\begin{align*}
i(j_{\bar F}(j_{\bar F})(\bar E)) &= (k\circ k_*)(j_{\bar F}(j_{\bar F})(\bar E))\\
&= k(k_*(j_{\bar F}(j_{\bar F})(\bar E)))\\
&= k(j_{\bar F}(k)(j_{\bar F}(j_{\bar F})(\bar E)))\\
&= k(j_{\bar F}(k)\circ (j_{\bar F}(j_{\bar F}))(\bar E))\\
&= k(j_{\bar F}(k\circ j_{\bar F})(\bar E))\\
&= k(j_{\bar F}(j_F)(\bar E))\\
&= k(j_{\bar F}(j_F))(k(\bar E))\\
&= j_F(j_F)(k(\bar E))\\
&= j_F(j_F)(E)
\end{align*}
This proves the claim.
\end{proof}
So the hypotheses of the theorem hold for \(\bar F\) and \(\bar E\). In other words, replacing \(F,E\) with \(\bar F,\bar E\), we may assume that \(F\) is an ultrafilter.

Assume towards a contradiction the theorem fails. Repeating the proof of \cref{rangefix}, there is an ordinal \(\xi\) and a sequence of extenders \(\langle E_n : n < \omega\rangle\) such that \(j_{E_n}(\xi) > \xi\) for all \(n\), and for all \(n < m\), \(E_n\) and \(E_m\) commute. But this is impossible by the proof of \cref{KunenClaim}.
\end{proof}
\end{thm}

A special case is the following corollary:

\begin{cor}
Suppose \(\mathcal U\) is a normal fine \(\kappa\)-complete ultrafilter on \(P_\kappa(\delta)\) and \(W\in \textnormal{Un}_\delta\). If \(W\mo \mathcal U\) then \[j_W\restriction \textnormal{Ord} \leq j_\mathcal U\restriction \textnormal{Ord}\] and in fact \(j_W\) fixes every ordinal in the range of \(j_\mathcal U\).\qed
\end{cor}

Also note that the corollary can {\it fail} if the commutativity hypothesis of \cref{rangefix2} fails. For example, if \(\kappa < \delta\) are measurable cardinals and \(U\) is a \(\delta\)-complete ultrafilter on \(\delta\) and \(\mathcal W\) is a \(\kappa\)-complete ultrafilter on \(P_\kappa(\delta)\) and \(\mathcal W \mo U\), then \(j_\mathcal W\restriction \text{Ord}\) is {\it not} dominated by \(j_U\restriction \text{Ord}\) since \(j_\mathcal W(\kappa) > \kappa= j_U(\kappa)\).

In fact, the converse of \cref{rangefix} is also true.

\begin{prp}[UA]
Suppose \(U\in \textnormal{Un}\), \(W\in \textnormal{Un}^{M_U}\), and  \[j^{M_U}_W\circ j_U\restriction \textnormal{Ord} = j_U\restriction \textnormal{Ord}\] Then in \(M_{U^2}\), \(j_U(W)\I j_U(j_U)(W)\).
\begin{proof}
Fix a sufficiently large cardinal \(\delta\). Let \(\xi = |j_U(j_U)(W)|^{M_{U^2}}_{\leq \delta}\). Then \(\xi = j_U(j_U)(|W|^{M_U}_{\leq\delta})\). Hence \(j_U(j^{M_U}_W)\) fixes \(\xi\), since by elementarity \(j_U(j^{M_U}_W)\) fixes every ordinal in the range of \(j_U(j_U)\). It follows from \cref{FixThm2} that that in \(M_{U^2}\), \(j_U(W)\I j_U(j_U)(W)\).
\end{proof}
\end{prp}

It is a bit bizarre that the question of whether \(U\I W\) depends only on the fixed points of \(j_U\). But perhaps assuming UA, there is some way to reconstruct an elementary embedding from its restrictions to large enough sets of ordinals. 
\begin{qst}[UA]
Suppose \(U\) and \(W\) are countably complete ultrafilters and there are arbitrarily large sets \(A\subseteq \textnormal{Ord}\) such that \(j_U\restriction A \in M_W\). Must \(U \I W\)?
\end{qst}

We now give a partial answer to this question:
\begin{thm}[UA]\label{OrdinalInternal}
Suppose \(U\) and \(W\) are countably complete ultrafilters with the property that \(j_U\restriction \alpha\in M_W\) for all ordinals \(\alpha\). Then \(U \I W\).
\end{thm}
This is not provable in ZFC. We sketch the independence result. Assume \(\kappa\) is measurable of Mitchell order 2 and fix normal ultrafilters \(U\mo W\) on \(\kappa\). By Kunen-Paris forcing \cite{KunenParis}, let \(N\supseteq V\) be a cardinal-preserving generic extension with the same continuum function such that every normal ultrafilter of \(V\) on \(\kappa\) lifts to \(2^{2^\kappa}\) normal ultrafilters in \(N\). Let \(W^*\in N\) be a lift of \(W\). Then \(M_W^V\subseteq M_{W^*}^N\) so \(j^V_U\restriction \alpha\in M_{W^*}^N\) for all ordinals \(\alpha\). Therefore for every lift \(U^*\) of \(U\) to \(N\), \(j_{U^*}^N\restriction \alpha = j^V_U\restriction \alpha\in M_{W^*}^N\). But since \((2^{2^\kappa})^{M_{W^*}^N} < (2^\kappa)^+\), not all such lifts \(U^*\) can belong to \(M_{W^*}^N\). 

For the proof we need the following fact:
\begin{lma}[UA]\label{MinftyLemma}
Suppose \(U\) and \(W\) are countably complete ultrafilters. Then there is an inner model \(Q\) admitting an elementary embedding \(j : V \to Q\) and internal elementary embeddings \((i_0,i_1) : (M_U,M_W)\to Q\) such that \(i_0\circ j_U = i_1 \circ j_W = j\) and such that \(j_U(Q) = j_W(Q) = Q\).
\end{lma}

To build \(Q\), we consider certain directed systems which are closely related to the proof of \cref{PointedAbs}:
\begin{defn}
If \(\lambda\) is a cardinal, then \(\mathcal D_\lambda\) denotes the category of ultrapowers of \(V\) by countably complete ultrafilters on \(\lambda\) with internal ultrapower embeddings.
\end{defn}

\begin{lma}[UA]\label{DAbs}
For any cardinal \(\lambda\), \(\mathcal D_\lambda\) is a directed partial order. For any ultrapower embedding \(j : V\to M\) with \(M\in \mathcal D_\lambda\), \((\mathcal D_{j(\lambda)})^M\) is equal to the collection of ultrapowers in \(\mathcal D_\lambda\) that lie above \(M\) in this partial order.
\begin{proof}
This is just like the argument for \cref{PointedAbs}.
\end{proof}
\end{lma}

Since \(\mathcal D_\lambda\) is a directed partial order, it makes sense to take the direct limit of the embeddings of \(\mathcal D_\lambda\).

\begin{defn}[UA]
For any cardinal \(\lambda\), let  \(M_\lambda = \lim \mathcal D_\lambda\). For any ultrapower \(N\in \mathcal D_\lambda\), let \(j_{N,\lambda}: N\to M_\lambda\) denote the direct limit embedding.
\end{defn}

It is worth mentioning the following fact though we will not use it here:
\begin{prp}[UA]
For any \(\mathcal M\in \mathscr P_\lambda\), \(|\mathcal M|_\lambda = j_{\mathcal M,\lambda}(\alpha_\mathcal M)\).\qed
\end{prp}

\begin{proof}[Sketch of \cref{MinftyLemma}]
We may assume without loss of generality that \(U\) and \(W\) are uniform ultrafilters on ordinals. Fix a cardinal \(\lambda\) such that \(U,W\in \Un_{\leq\lambda}\). Then \(j_U(M_\lambda) = M_{j_U(\lambda)}^{M_U} = \lim (\mathcal D_{j_U(\lambda)})^{M_U}\). But by \cref{DAbs}, \((\mathcal D_{j_U(\lambda)})^{M_U}\) is the cone above \(M_U\) in the directed system \(\mathcal D_\lambda\). Therefore \( \lim  (\mathcal D_{j_U(\lambda)})^{M_U} = \lim \mathcal D_\lambda = M_\lambda\). This shows \(j_U(M_\lambda) = M_\lambda\). Similarly \(j_W(M_\lambda) = M_\lambda\). Letting \(Q = M_\lambda\), \(i_0 = j_{M_U,\lambda}\), and \(i_1 = j_{M_W,\lambda}\), this proves the theorem.
\end{proof}

\begin{proof}[Proof of \cref{OrdinalInternal}]
Let \(j : V \to Q\) be as in \cref{MinftyLemma} and fix an elementary embedding \(i : M_W\to Q\) that is definable over \(M_W\) and satisfies \(i\circ j_W = j\). We claim \(j_U\restriction Q\) is amenable to \(M_W\) in the sense that for any \(X\in Q\), \(j_U\restriction X\in M_W\). To see this, let \(f : \alpha\to X\) be a bijection between an ordinal \(\alpha\) and \(X\) such that \(f\in Q\). Let \(X' = j_U(X)\) and let \(f' = j_U(f)\). Since \(j_U(Q)\subseteq Q\), \(X'\) and \(f'\) belong to \(Q\). But \[j_U\restriction X = f'\circ (j_U\restriction \alpha)\circ f^{-1}\] and so \(j_U\restriction X\in M_W\) since the functions on the righthand side belong to \(M_W\).

Now let \[U_* = \{X\in P^Q(\sup j[\delta]) : j^{-1}[X]\in U\}\] By the proof of \cref{Pushforward}, \(U_*\) is the uniform \(Q\)-ultrafilter derived from \(j_U\restriction Q\) using \(j_U(j)([\text{id}]_U)\). Therefore \(U_*\in M_W\) since \(j_U\restriction Q\) is amenable to \(M_W\). But note that 
\begin{align*}X\in \s W U&\iff j_W^{-1}[X]\in U\\
&\iff (i\circ j_W)^{-1}[i(X)]\in U\\
&\iff j^{-1}[i(X)]\in U\\
&\iff i(X)\cap \sup j[\delta]\in U_*
\end{align*}
Since \(i\) is definable over \(M_W\) and \(U_*\in M_W\), it follows that \(\s W U\) can be computed inside \(M_W\). Therefore by \cref{Pushforward}, \(U\I W\).
\end{proof}
\section{Questions}
We pose two questions related to the internal relation. Our first is whether \cref{CommutingConverse} is provable in ZFC:
\begin{qst}
Suppose \(U\) and \(W\) are countably complete ultrafilters such that \(U\I W\) and \(W\I U\). Must \(j_U(j_W) = j_W\restriction M_U\) and \(j_W(j_U) = j_U\restriction M_W\)?
\end{qst}

Our last question is related to \cref{SelfAmenable}:
\begin{qst}
What is the consistency strength of the existence of a nonprincipal countably complete ultrafilter on a cardinal that is amenable to its own ultrapower? Can such an ultrafilter be weakly normal? What if UA holds?
\end{qst}

\bibliography{Bibliography}{}
\bibliographystyle{unsrt}

\end{document}